%% file: test.tex
\documentclass[11pt]{article}
\usepackage{epic}
\usepackage{html}

\include{defines}

\begin{document}
\title{A non-regular Gr\"obner fan}
\author{Anders Nedergaard Jensen
\thanks{Research partially supported by the Faculty of Science, University of Aarhus, Danish Research Training Council (Forskeruddannelsesr\aa det, FUR) , Institute for Operations Research ETH, grants DMS 0222452 and DMS 0100141 of the U.S. National Science Foundation and the American Institute of Mathematics.
}
\\
\\
\small
Department of Mathematical Sciences, University of Aarhus
and\\
\small
Institute for Operations Research, ETH Z\"urich
}

\maketitle
\begin{abstract}
  The Gr\"obner fan of an ideal $I\subset k[x_1,\dots,x_n]$, defined
  by Mora and Robbiano, is a complex of polyhedral cones in $\R^n$.
  The maximal cones of the fan are in bijection with the distinct
  monomial initial ideals of $I$ as the term order varies. If $I$ is
  homogeneous the Gr\"obner fan is complete and is the normal fan of
  the state polytope of $I$.
  In general the Gr\"obner fan is not complete and therefore not the
  normal fan of a polytope. We may ask if the \emph{restricted}
  Gr\"obner fan, a subdivision of $\R_{\geq 0}^n$, is regular i.e. the
  normal fan of a polyhedron. The main result of this paper is an
  example of an ideal in $\Q[x_1,\dots,x_4]$ whose restricted
  Gr\"obner fan is not regular.
\end{abstract}

\section{Introduction}
Let $R=k[x_1,\dots,x_n]$ be the polynomial ring in $n$ variables over
a field $k$ and let $I\subset R$ be an ideal. The \emph{Gr\"obner fan} and the \emph{restricted Gr\"obner fan} of $I$ are $n$-dimensional polyhedral fans defined in \cite{MoRo}. 
The main result of this paper is the following.

\begin{theorem}
\label{main theorem}
The restricted Gr\"obner fan of the two-dimensional ideal
$$I=\langle acd+a^2c-ab, ad^2-c, ad^4+ac\rangle\subset\Q[a,b,c,d]$$
is not the normal fan of a polyhedron.
\end{theorem}

In contrast, when the ideal $I$ is \emph{homogeneous} its Gr\"obner fan and restricted Gr\"obner fan are known to be normal fans of polyhedra, see Section \ref{section definition}.

We recall the definition of a \emph{fan} in $\R^n$. A \emph{polyhedron} in $\R^n$ is a set of the form $\{x\in\R^n:Ax\leq b\}$ where $A$ is a matrix and $b$ is a vector. Bounded polyhedra are called \emph{polytopes}. If $b=0$ the set is a \emph{polyhedral cone}.

\begin{definition}
A collection $C$ of polyhedra in $\R^n$ is a polyhedral complex if:
\begin{enumerate}
\item all proper faces of a polyhedron $P\in C$ are in $C$, and
\item the intersection of any two polyhedra $A,B\in C$ is a face of $A$ and a face of $B$. 
\end{enumerate}
\end{definition}
A polyhedral complex is a \emph{fan} if it only consists of cones. A simple way to construct a fan is by taking the \emph{normal fan} of a polyhedron.
\begin{definition}
Let $P\subset\R^n$ be a polyhedron. All non-empty faces of $P$ are of the form
$$ face_\omega(P)=\{p\in P:\langle \omega,p\rangle=\textup{max}_{q\in P}\langle \omega,q\rangle\}$$
for some $\omega\in\R^n$.
For a face $F$ of $P$ we define its normal cone
$$N_P(F):=\overline{\{\omega\in\R^n:face_\omega(P)=F\}}$$
with the closure being taken in the usual topology. The normal fan of $P$ is the fan consisting of the normal cones $N_P(F)$ as $F$ runs through all non-empty faces of $P$.
\end{definition}
If the union of all cones in a fan is $\R^n$, the fan is said to be \emph{complete}. It is clear that the normal fan of a polytope is complete.
Not all fans arise as the normal fan of a polyhedron. Those that do are called \emph{regular}.

If the ideal $I$ is homogeneous, its Gr\"obner fan is the normal fan of a
polytope known as the \emph{state polytope} of $I$ (\cite{bayer}, \cite[Chapter 2]{sturmfels}). In the general
case, no similar result exists as the Gr\"obner fan is not
complete.
However, we could ask if the \emph{restricted} 
Gr\"obner fan of $I$, a fan in
$\R_{\ge 0}^n$, is regular. Theorem \ref{main
theorem} gives an example of an ideal in $\Q[x_1,\dots,x_4]$ whose
Gr\"obner fan and restricted Gr\"obner fan are not regular.

The definitions of the Gr\"obner fan and the restricted Gr\"obner fan appear in Section \ref{section definition}, and the proof of Theorem \ref{main theorem} is given in Section \ref{section example}. For the reader unfamiliar with Gr\"obner fans we provide the necessary background in Section \ref{section definition}.
It is interesting to consider what happens if we homogenize the example ideal $I$ and project its state polytope back into $\R^4$. In Section \ref{section comments} we will point out why the normal fan of this projection is not the Gr\"obner fan of $I$. In particular we conclude for this example that the third variant of the Gr\"obner fan, the \emph{extended} Gr\"obner fan defined in \cite{MoRo}, does not agree with the restricted fan in the positive orthant.

An interesting corollary of the restricted Gr\"obner fan being regular would be an easy proof that the memoryless reverse search algorithm (\cite{avis}) can be used for enumerating the maximal cones in the fan by exploiting the structure of the underlying polyhedron. In light of Theorem \ref{main theorem} the fact that the reverse search method can be used requires a non-trivial proof which will appear in \cite{fukuda}.

\section{The Gr\"obner fan of an ideal}
\label{section definition}
For $\alpha\in\N^n$ we use the notation $x^\alpha:=x_1^{\alpha_1}\dots x_n^{\alpha_n}$. By
a term order on $R$ we mean a total ordering on the monomials in $R$
such that:
\begin{enumerate}
\item For all $\alpha\in\N^n\backslash \{0\}:1<x^\alpha$ and
\item For $\alpha,\beta,\gamma\in\N^n: x^\alpha < x^\beta \Rightarrow x^\alpha x^\gamma < x^\beta x^\gamma$.
\end{enumerate}
Let $\prec$ be a term order. For a non-zero polynomial $f\in R$ we
define its \emph{initial term}, $in_\prec(f)$, to be the unique maximal
term of $f$ with respect to $\prec$. In the same way for $\omega\in
\R^n$ we define the \emph{initial form}, $in_\omega(f)$, to be the sum
of all terms of $f$ whose exponents maximize $\langle \cdot,
\omega\rangle$. The \emph{initial ideals} of an ideal $I$ with respect to $\prec$ and $\omega$ are defined
as
$$ in_\prec(I)=\langle in_\prec(f):f\in I\backslash\{0\}\rangle \,\,\mbox{and}\,\, in_\omega(I)=\langle in_\omega(f):f\in I\backslash\{0\}\rangle.$$
Note that $in_\prec(I)$ is a monomial ideal while $in_\omega(I)$ might not be. A monomial in $R\backslash in_\prec(I)$ (with coefficient $1$) is called a \emph{standard monomial} of $in_\prec(I)$.
\begin{definition}
Let $I\subset R$ be an ideal and $\prec$ a term order on $R$. A generating set $\G=\{g_1,\dots,g_m\}$ for $I$ is called a Gr\"obner basis for $I$ with respect to $\prec$ if
$$in_\prec(I)=\langle in_\prec(g_1),\dots,in_\prec(g_m)\rangle.$$
The Gr\"obner basis $\G$ is minimal if no polynomial can be left out. A minimal Gr\"obner basis is reduced if the initial term of every $g\in\G$ has coefficient $1$ and all other monomials in $g$ are standard monomials of $in_\prec(I)$.
\end{definition}
For a term order $\prec$ and an ideal $I$ the reduced Gr\"obner basis is unique and depends only on $I$ and $in_\prec(I)$. We denote it by $\G_\prec(I)$.

Given $I$ a natural
equivalence relation on $\R^n$ is the one induced by taking initial
ideals:
$$u\sim v \Leftrightarrow in_u(I)=in_v(I).$$
We introduce the following notation for the closures of the equivalence classes:
$$C_\prec(I)=\overline{\{u\in\R^n : in_u(I)=in_\prec(I)\}}\mbox{~~and}$$
$$C_v(I)=\overline{\{u\in\R^n : in_u(I)=in_v(I)\}}.~~~~~$$
A well known fact is that for a
fixed ideal $I$ there are only finitely many sets $C_\prec(I)$  and they cover $\R_{\geq 0}^n$, see \cite{MoRo}. Secondly,
every initial ideal $in_\prec(I)$ is of the form $in_\omega(I)$
for some $\omega\in\R_{>0}^n$. Consequently, every $C_\prec(I)$ is of the form
$C_\omega(I)$.
A third observation is that the equivalence classes are not convex in
general since we allow the vectors to be anywhere in $\R^n$:
\begin{example}
\label{first example}
Let $I=\langle x_1-1, x_2-1\rangle$. The ideal $I$ has five initial
ideals: $\langle x_1-1, x_2-1\rangle$, $\langle x_1, x_2\rangle$,
$\langle x_1, x_2-1\rangle$, $\langle x_1-1, x_2\rangle$ and $\langle
1\rangle$. In particular, for $u=(-1,3)^T$ and $v=(3,-1)^T$ we have
$in_u(I)=in_v(I)=\langle 1\rangle$ but $in_{{1\over
2}(u+v)}(I)=\langle x_1, x_2 \rangle$.
\end{example}
\begin{theorem}
\label{groebner cone theorem}
Let $\prec$ be a term order and $v\in C_\prec(I)$ then for
$u\in\R^n$
$$in_u(I)=in_v(I) ~\Longleftrightarrow ~\forall g\in\G_\prec(I):in_u(g)=in_v(g).$$
\end{theorem}
This theorem is a little more general than Proposition 2.3 in
\cite{sturmfels} as it allows the vectors to be negative. A proof will
appear in \cite{fukuda}. Theorem \ref{groebner cone theorem} shows that the closures of the
equivalence classes are polyhedral cones since for fixed $\prec$ and
fixed $v$ each $g\in\G_\prec(I)$ introduces the equality $in_u(g)=in_v(g)$
which is equivalent to having $u$ satisfy a set of linear equations
and strict linear inequalities. The closure is taken by making the
strict inequalities non-strict.
Thus in particular, the set $C_v(I)$ is a convex polyhedral cone if it contains a strictly positive vector.

\begin{definition}
\label{definition gfan}
The \emph{Gr\"obner fan} of an ideal $I\subset R$ is the set of the
closures of all equivalence classes intersecting the positive orthant
together with their proper faces.
\end{definition}
This is a variation of the definitions appearing in the
literature. The advantage of this variant is that it gives
well-defined and \emph{nice} fans in the homogeneous and
non-homogeneous case simultaneously. By \emph{nice} we mean that all cones in this fan are closures of equivalence classes.
It is not clear a priori that the Gr\"obner fan is a polyhedral complex. The proof that the Gr\"obner fan is in fact a fan (polyhedral complex) will be deferred to \cite{fukuda}. The support of the Gr\"obner fan of $I$ is called the \emph{Gr\"obner region} of $I$.

For the purpose of this paper it is better to study the
\emph{restricted} Gr\"obner fan as we will see soon. Using the
definition we already have together with the notion of \emph{common
refinements} of fans (\cite{ziegler}) it is straightforward to make a
definition equivalent to the original one in \cite{MoRo}.
\begin{definition}
Let $F$ and $F'$ be two polyhedral fans in $\R^n$. Their common refinement is the polyhedral fan $F\wedge F':=\{C\cap C'\}_{(C,C')\in F\times F'}$.
\end{definition}
\begin{definition}[Definition 2.5 \cite{MoRo}]
\label{definition rgfan}
The restricted Gr\"obner fan of an ideal $I\subset R$
is the common refinement of the non-negative orthant $\R_{\ge 0}^n$
with its proper faces and the Gr\"obner fan of $I$.
\end{definition}
The support of the restricted Gr\"obner fan is $\R_{\geq 0}^n$.

A fundamental question to ask is the following: \emph{Is the Gr\"obner
fan always the normal fan of a polytope?} The answer to this question is \emph{no} since the Gr\"obner fan is not always complete.
Even if we ask for a polyhedron instead, the answer is still \emph{no}
for trivial reasons as the following example shows.
 \begin{example}
Let $f=x_1+x_2+1$ and $I=\langle f\rangle\subset\Q[x_1,x_2]$.
The picture to the left shows the two maximal cones in the Gr\"obner
fan of $I$. Any polyhedron having these cones as normal cones will have at least a third normal cone (middle picture).
\begin{center}
\begin{minipage}[0cm]{12.6cm}
\epsfig{file=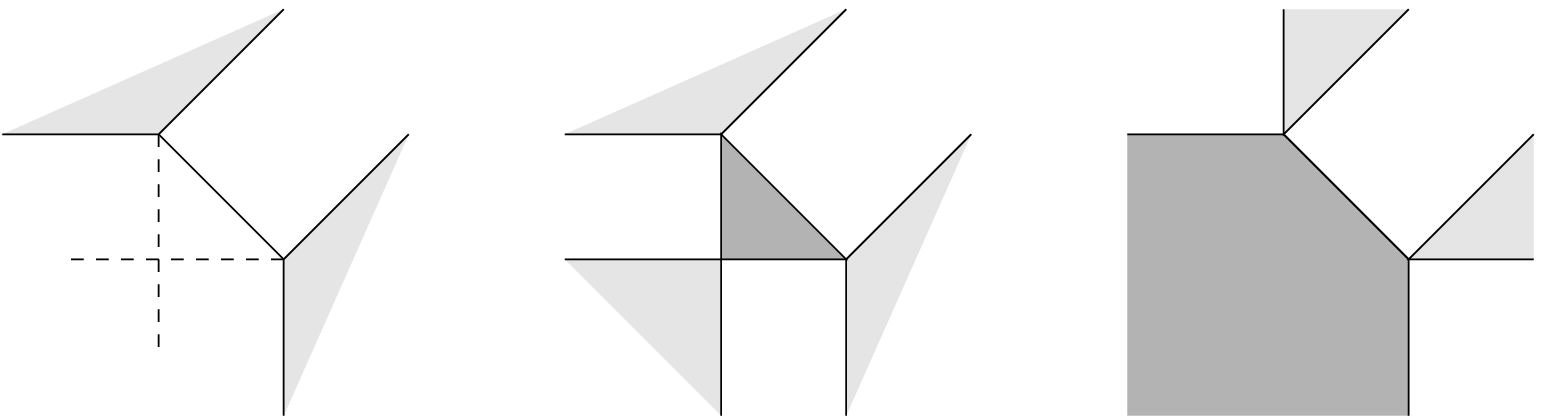,height=3.4cm} 
\end{minipage}
\end{center}
To the right the restricted Gr\"obner fan is shown. In this example it is the normal fan of an unbounded polyhedron.
\end{example}

Thus we rephrase the question for restricted Gr\"obner fans:
\emph{Is the restricted Gr\"obner fan of an ideal always the normal
fan of a polyhedron?}

We note that the Gr\"obner fan being regular is
stronger than the restricted Gr\"obner fan being so. This is because
the normal fan of the Minkowski sum of two polyhedra is the
common refinement of their normal fans. The claim follows since
$\R_{\geq 0}^n$ with its proper faces is the normal fan of $\R_{\leq
0}^n$.

The above question is known to have a positive answer in the following three special cases:
\begin{itemize}
\item 
If the ideal is homogeneous the answer is \emph{yes} since the
Gr\"obner fan is the normal fan of the \emph{state polytope} of $I$ introduced by Bayer and Morrison in \cite{bayer}.
We should mention that in \cite{mall} it is shown that the Gr\"obner fan
is not the normal fan of the state polytope
as it was defined in
\cite{bayer}. Instead we should use the construction
in \cite[Chapter 2]{sturmfels}.  We take the Minkowski sum
of the state polytope with $\R_{\leq 0}^n$ to get a polyhedron having the restricted
Gr\"obner fan as its normal fan.

\item
The \emph{Newton polytope}, $New(f)$, of a polynomial $f$ is defined to
be the convex hull of the exponent vectors of the monomials in $f$. In the
case of a \emph{principal} ideal $I=\langle f\rangle$ the Newton
polytope $New(f)$ will almost have the Gr\"obner fan as its normal fan
since two vectors $u,v\in\R^n$ pick out the same initial ideal of $I$
if and only if they are maximized on the same face of $New(f)$. The
only thing that keeps $New(f)$ from having the Gr\"obner fan of $I$ as
its normal fan is that we have not included all equivalence classes in
the Gr\"obner fan.  However, the normal fan of the Minkowski sum of
$New(f)$ and $\R_{\leq 0}^n$ is the restricted Gr\"obner fan.

\item
A third case where we have a similar result is for
\emph{zero-dimensional} ideals. The construction of a polytope is
similar but simpler than the construction in the homogeneous case as
there are only a finite number of standard monomials for each initial ideal. We
claim, without proof, that the following construction works: For every
term order $\prec$ construct the vector $v_\prec$ equal to the negative of the sum of all exponent
vectors of all standard monomials of $in_\prec(I)$. Take the convex hull of all
$v_\prec$ as we vary the term order. The Minkowski sum of this polytope with $\R_{\geq 0}^n$ is a polyhedron whose normal
fan is the restricted Gr\"obner fan.
\end{itemize}
In contrast to the above, we have Theorem \ref{main theorem}.

\section{The proof}
\label{section example}
This section contains a proof of Theorem \ref{main theorem}. We start
by deducing a necessary condition for a fan to be the normal fan of a polyhedron. We then
show that the restricted Gr\"obner fan of the ideal in the theorem violates
this condition.
 Finally we argue that the Gr\"obner fan has been
computed correctly.

\subsection{A necessary condition}
Let $F$ be a fan in $\R^n$. Suppose $F$ is the normal fan of a
polyhedron $P\subset\R^n$. The non-empty faces of $P$ are in bijection
with the cones in $F$ by taking normal cones of the faces. Adjacency
is preserved in the sense that two vertices of an edge of $P$ map to cones
in $F$ having the normal fan of the edge as a common facet. Furthermore, the
edge is perpendicular to the shared facet. If a set of normals of the
shared facets in $F$ are specified, then for every bounded edge the
difference between its endpoints can be expressed as some scalar
times the specified normal of its normal cone. The scalars are considered to be unknowns. Since the adjacency information of the vertices of $P$ is present in
$F$, the bounded edge graph of $P$ can be deduced from $F$. A necessary
condition for $F$ to be the normal fan of $P$ is that every
combinatorial cycle in the edge graph is a geometric cycle in
space. This condition gives rise to a feasible system of inequalities
on the scalars dependent on $F$ alone.

To be more specific about the inequality system, consider the adjacency
graph of the $n$-dimensional cones in $F$, or equivalently the edge
graph of the supposed polyhedron $P$. Let $V=\{1,\dots,m\}$ denote the
vertices and a subset $E\subset\{(i,j)\in V\times V:i<j\}$ denote the
edges in the graph. For each shared facet, choose a normal vector $d_{(i,j)}\in\R^n$ such that the $i$th cone is
on the negative side of the hyperplane with normal vector $d_{(i,j)}$ and
the $j$th cone is on the positive side. The graph $(V,E)$ is considered to be
undirected when we define its cycles. A vector $f\in\R^E$ is called a
flow in $(V,E)$ if
$$\forall j\in V:\sum_{(i,j)\in E}f_{(i,j)}=\sum_{(j,k)\in
  E}f_{(j,k)}.$$
In other words the flow entering $j$ is the same as the
flow leaving $j$. The set of flows is a subspace of $\R^E$.  We
introduce a vector $s\in\R_{>0}^E$ of unknown scalars such that the true vector
from vertex $i$ to vertex $j$ is $s_{(i,j)}d_{(i,j)}$. Each cycle in
the graph can be represented by a flow $f\in\R^E$ being $0$ on the edges
not appearing in the cycle and $\pm 1$ elsewhere depending on the
relative orientation of the cycle and the edge. For such an $f$ the
condition that the cycle forms a loop in space can be expressed as:
\begin{equation}
\sum_{(i,j)\in E}f_{(i,j)}s_{(i,j)}{ d}_{(i,j)}={ 0} \label{system}.
\end{equation}
Note that (\ref{system}) is a system of $n$ equations -- one for each coordinate of $d_{(i,j)}$. If $F$ is the normal fan of
a polyhedron $P$, there exist positive scalars $s_{(i,j)}$ satisfying (\ref{system}) for
every flow $f$ since the cycle flows span the vector space of flows.  By
linearity this is equivalent to having the scalars satisfy (\ref{system}) for
a basis of the vector space of flows rather than the entire space. In
matrix form we may express the necessary condition as the system
\begin{equation}
\label{matrix system}
As=0  \mbox{~ and ~} s_{(i,j)}>0 ~\mbox{for all}~ (i,j)\in E
\end{equation}
having a solution where $A$ is a suitable $nl\times|E|$ matrix with
$l$ being the dimension of the vector space of flows.

\subsection{The certificate}
\begin{figure}
\vspace{1.1cm}
\noindent
\hspace{0.15cm}
\input{testfig}
\caption{The certificate subgraph.\label{the figure}}
\end{figure}
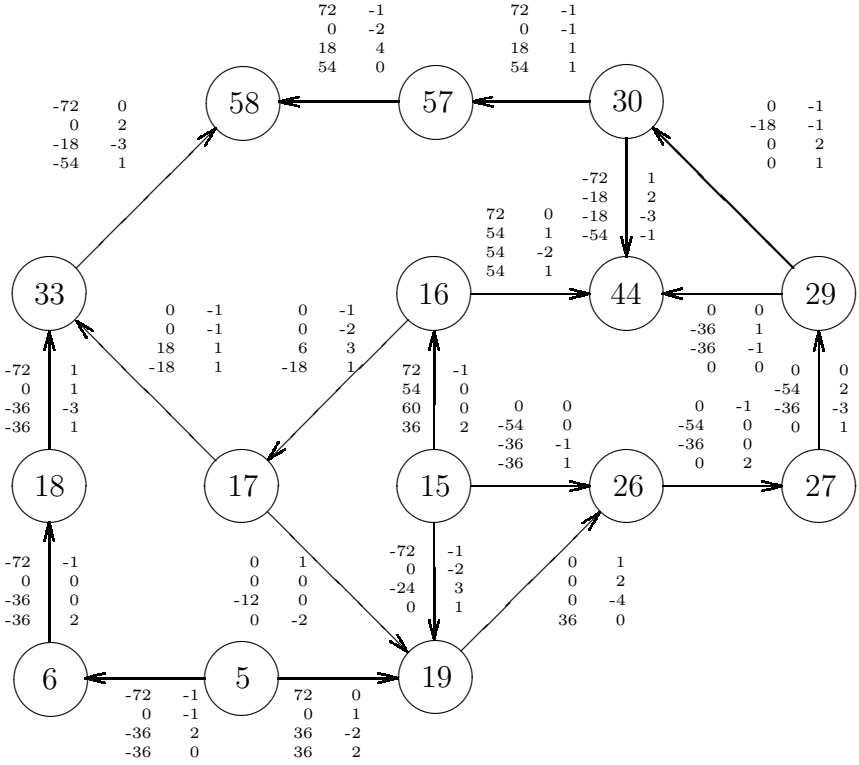

\noindent
{\it Proof of Theorem \ref{main theorem}:}
\label{the example}
The restricted Gr\"obner fan of the ideal
$$I=\langle acd+a^2c-ab, ad^2-c, ad^4+ac\rangle\subset\Q[a,b,c,d]$$
has 81 full dimensional cones each corresponding to a monomial initial
ideal. Their adjacency graph $(V,E)$ has 163 edges, with each edge having an
edge direction equal to a specified normal of the shared facet. We
present a certificate that the fan is not the normal fan of a
polyhedron.
Only the subgraph in Figure \ref{the figure} is needed to describe it.
Two vectors are written for each edge in the subgraph. The vector to
the right is the edge direction $d_{(i,j)}$ and the vectors to the left
describe four flows in the subgraph.

Let $V'$ be the set of vertices appearing in the subgraph and $E'$ the
edges. Let $f^1,f^2,f^3$ and $f^4$ denote the flows above. Suppose the
restricted Gr\"obner fan was the normal fan of a polyhedron
$P$. Equality system (\ref{system}) implies
\begin{equation}
\label{system2}
\forall (r,t)\in\{1,2,3,4\}\times\{1,2,3,4\}:\sum_{(i,j)\in E'}f^r_{(i,j)}s_{(i,j)}{ d}_{(i,j)_t}=0.
\end{equation}
In particular, the sum of the equations in (\ref{system2}) for $(r,t)=(1,1),(2,2),(3,3),(4,4)$ is zero. Therefore,
$$0=\sum_{r=1}^4\sum_{(i,j)\in E'} s_{(i,j)} d_{{(i,j)}_r} f_{{(i,j)}}^r=\sum_{{(i,j)}\in E'}s_{(i,j)}\sum_{r=1}^4  d_{{(i,j)}_r} f_{{(i,j)}}^r.
$$
The local contribution at each edge except the edge (29,30) is zero
because $d_{(i,j)}\cdot
(f_{(i,j)}^1,f_{(i,j)}^2,f_{(i,j)}^3,f_{(i,j)}^4)^T=0$ (check this in
the picture). Consequently,
$$0=s_{(29,30)}d_{(29,30)}\cdot f_{(29,30)}=18s_{(29,30)}$$
implying $s_{(29,30)}=0$. Hence the vertices 29 and 30 have the same coordinates which contradicts that $P$ is a polyhedron with the required edge graph.
$_\Box$

\begin{remark}
Another way to argue is by observing that we have applied the trivial
direction of Farkas' lemma to (\ref{system2}). With
$A'$ being the $16 \times 20$ matrix representing the equalities in
(\ref{system2}) a variant of Farkas' lemma says:
$$\exists ~y : y^T A' \ge 0 \mbox{ and } y^T A' \neq 0
~\Longleftrightarrow ~ \not\exists ~s>0: A's=0 .$$ In our case
$y=(1,0,0,0,0,1,0,0,0,0,1,0,0,0,0,1)^T$ where the four nonzero
components correspond to the equations $(1,1)$, $(2,2)$, $(3,3)$
and $(4,4)$.
\end{remark}

\subsection{Correctness of the subgraph}
For completeness, a positive interior point in each of the 15 maximal
cones of the restricted Gr\"obner fan leading to the inconsistency is given in the top part of Figure \ref{the second figure}. Further, a positive
vector in the relative interior of every shared facet is given in the bottom part.

\begin{figure}
\begin{center}
\begin{tabular}{|c|c|c|}
\hline
\begin{tabular}{rl}
5 & $(10,2,5,3)$
 \\
6 & $(14,4,11,5)$
 \\
15 & $(7,6,5,3)$
 \\
16 & $(7,11,8,4)$
 \\
17 & $(5,2,3,3)$
 \\
\end{tabular}
&
\begin{tabular}{rl}
18 & $(4,3,5,4)$
 \\
19 & $(5,1,2,2)$ 
 \\
26 & $(7,1,2,3)$
 \\
27 & $(17,1,4,9)$
 \\
29 & $(10,1,2,6)$
 \\
\end{tabular}
&
\begin{tabular}{rl}
30 & $(15,1,3,11)$
 \\
33 & $(3,1,2,3)$
 \\
44 & $(7,5,4,4)$
 \\
57 & $( 7,1,2,7)$ 
 \\
58 & $(7,1,3,8)$ 
 \\
\end{tabular}\\
\hline
\begin{tabular}{rrl}
5 & 6  & (3,1,2,1)\\
5 & 19 & (8,4,5,3)\\
6 & 18 & (2,1,2,1)\\
15 & 16 & (6,8,6,3)\\
15 & 19 & (5,3,3,2)\\
15 & 26 & (9,2,3,3)\\
16 & 17 & (8,15,11,5)\\
\end{tabular}
&
\begin{tabular}{rrl}
16 & 44 & (5,7,5,3)\\
17 & 19 & (4,1,2,2)\\
17 & 33 & (6,1,3,4)\\
18 & 33 & (4,1,3,4)\\
19 & 26 & (10,1,3,4)\\
26 & 27 & (18,2,5,9)\\
27 & 29 & (13,1,3,7)\\
\end{tabular}
&
\begin{tabular}{rrl}
29 & 30 & (8,1,2,5)\\
29 & 44 & (9,3,3,5)\\
30 & 44 & (6,5,4,4)\\
30 & 57 & (13,1,3,11)\\
33 & 58 & (6,1,3,7)\\
57 & 58 & (10,1,3,11)\\
& &\\
\end{tabular}
\\
\hline
\end{tabular}
\end{center}
\caption{Representative weight vectors for cones in the certificate.\label{the second figure}}
\end{figure}
To verify the correctness of the certificate the following procedure
is suggested: It is straightforward to check that the flows are flows
and that the dot products of flows and listed directions are 0 except
for the edge (29,30). The question is how to check the correctness of the edge
subgraph and the listed directions. For each of the listed edges
$(i,j)$ with $i<j$ compute the corresponding reduced Gr\"obner bases
$\G_i$ and $\G_j$ and use Theorem \ref{groebner cone theorem} to
compute their cones $C_i$ and $C_j$. Check that the listed facet
vector for the edge $(i,j)$ is in the closure of both cones $C_i$ and $C_j$
and that the listed direction vector non-strictly separates $C_i$ and
$C_j$ with $C_j$ being on the non-negative side.  Checking that the
listed facet vector is in the relative interior of a facet of $C_r$
completes the verification. The non-straightforward part of this test
was implemented as a 230 line script in Singular \cite{singular}. The
script itself is available on the internet, see \cite{singularscript}.

\section{Further remarks}
\label{section comments}
\subsection{Homogenizing the ideal}
In \cite{MoRo} a complete fan in $\R^{n}$ called
the \emph{extended} Gr\"obner fan is defined for any (not necessarily homogeneous) ideal $I\subset
R$. This is done by homogenizing the ideal with a new variable. The \emph{extended} Gr\"obner fan
is defined as the Gr\"obner fan of the homogenized ideal intersected
with $\R^n$. It is clear that the extended Gr\"obner fan is regular
as the Gr\"obner fan of the homogenized ideal is regular and the normal fan of the projection of its polytope to $\R^n\times\{0\}$
is the intersection of the Gr\"obner fan of the homogenized ideal with $\R^n\times\{0\}$.
Therefore our example shows that the restricted Gr\"obner fan of an ideal and its
extended Gr\"obner fan need not agree in $\R_{>0}^n$.

In our example the procedure works as follows. We homogenize the ideal $I$ using the variable ``e'' to get
$$
^hI=\langle cd^{2}+ace,
-c^{2}e+c^{2}d+abd,
c^{2}e+c^{3}-bce-bcd-abd+abc,
-ce^{2}+ad^{2},$$
$$-c^{2}e+acd-abe,
c^{2}e-bce+ac^{2}-abd,
c^{2}e+a^{2}c,
bce+a^{2}b\rangle.
$$
The Gr\"obner fan of the new ideal is a complete fan in $\R^5$. Intersecting this fan with $\R^4\times\{0\}$ we get the extended Gr\"obner fan, a regular fan that almost equals the Gr\"obner fan of $I$ in the positive orthant. The subgraph listed for $I$ is valid for the extended fan on all edges except the edge connecting vertex $57$ and vertex $58$. The vector $(10,1,3,11)^T$ listed as a relative interior facet vector in the Gr\"obner fan of $I$ is not in the boundary of the cone containing $(7,1,2,7)^T$ in the extended fan.
 
\subsection{A program for finding the example}
A C++ program was written for finding non-regular Gr\"obner fans.
The input for the program is a set of generators for an
ideal $I$ and the output is either a coordinatization of a polyhedron
with the restricted Gr\"obner fan as its normal fan
or a certificate for its non-existence. The program works in two steps.
\begin{itemize}
\item In step 1 it calls a software package being developed by the author for computing Gr\"obner fans of polynomial ideals. This work will appear in \cite{fukuda}. The package computes the maximal cones ($n$-dimensional) of the Gr\"obner
fan of $I$ storing all facets ($n-1$-dimensional). This is done using
exhaustive search on the graph whose vertices are the maximal cones of
the fan, with two maximal cones being connected if they share a
facet. At each maximal cone the reduced Gr\"obner basis is known, its
facets are computed using linear programming and the
Gr\"obner bases of its neighbors are computed using the local basis
change procedure in \cite{collart}. A specialized implementation for
toric ideals was worked out in \cite{huber}.
\item From the Gr\"obner fan computed above the inequality system (\ref{matrix system})
is deduced. Linear programming methods are used for checking its
feasibility. The result is either positive scalars leading to a
coordinatization of the vertices of the polyhedron or a certificate for
its non-existence.
\end{itemize}
The software libraries \cite{gmp} and \cite{cdd} were used for doing the arithmetic and solving linear programming problems, respectively.

Knowing that we should avoid homogeneous, zero-dimensional and
principal ideals, it was not hard to find the example when the
C++ program had been written. A practical issue is that we
are restricted to ideals with not too complex Gr\"obner fans as the
entire edge graph must be handled by the LP-code. In looking for a
3-variable example this seems to be an unfortunate restriction as
nothing interesting happens in the small manageable examples we have
tried.

\vspace{1cm}
\noindent
{\bf Acknowledgments:}
The author is thankful to the following people and institutions for supporting this research: Komei Fukuda and Hans-Jakob L\"uthi (Institute for Operations Research, ETH Z\"urich), Douglas Lind and Rekha Thomas (University of Washington, Seattle) and the American Institute of Mathematics. In the writing process of this paper Niels Lauritzen, Komei Fukuda and, especially, Rekha Thomas have been very helpful. Thanks also to the many people who proofread this paper.

\bibliographystyle {plain}
\bibliography{jensen.bib}
 
\end{document}

%% file: defines.tex
\usepackage{theorem}
\usepackage[english]{babel}
\usepackage{a4}
\usepackage{epsfig}
\usepackage{amsfonts}
\usepackage{latexsym}
\newtheorem{theorem}{Theorem}

\newtheorem{definition}[theorem]{Definition}
\theorembodyfont{\rmfamily}
\newtheorem{remark}[theorem]{Remark}
\theorembodyfont{\rmfamily}
\newtheorem{example}[theorem]{Example}

\def\N{\ensuremath{\mathbb{N}}}

\def\Q{\ensuremath{\mathbb{Q}}}
\def\R{\ensuremath{\mathbb{R}}}

\def\G{\ensuremath{\mathcal{G}}}

%% file: testfig.tex
\setlength{\unitlength}{0.00083333in}
\begingroup\makeatletter\ifx\SetFigFont\undefined%
\gdef\SetFigFont#1#2#3#4#5{%
  \reset@font\fontsize{#1}{#2pt}%
  \fontfamily{#3}\fontseries{#4}\fontshape{#5}%
  \selectfont}%
\fi\endgroup%
{\renewcommand{\dashlinestretch}{30}
\begin{picture}(5408,4405)(0,-10)
\put(2775,2956){\circle{450}}
\put(2775,1756){\circle{450}}
\put(3975,1756){\circle{450}}
\put(3975,2956){\circle{450}}
\put(5175,1756){\circle{450}}
\put(5175,2956){\circle{450}}
\put(2784,565){\circle{450}}
\put(1575,556){\circle{450}}
\put(1575,1756){\circle{450}}
\put(384,547){\circle{450}}
\put(375,1756){\circle{450}}
\put(375,2956){\circle{450}}
\put(1584,4147){\circle{450}}
\put(2787,4152){\circle{450}}
\put(3975,4156){\circle{450}}
\drawline(2805.838,2602.960)(2776.000,2723.000)(2745.838,2603.041)
\drawline(2776,2723)(2775,1981)
\drawline(3000,1756)(3745,1755)
\drawline(3624.960,1725.161)(3745.000,1755.000)(3625.040,1785.161)
\drawline(4200,1756)(4945,1755)
\drawline(4824.960,1725.161)(4945.000,1755.000)(4825.040,1785.161)
\drawline(5175,1981)(5176,2721)
\drawline(5205.838,2600.960)(5176.000,2721.000)(5145.838,2601.041)
\drawline(3000,2956)(3745,2955)
\drawline(3624.960,2925.161)(3745.000,2955.000)(3625.040,2985.161)
\drawline(4320.000,2986.000)(4200.000,2956.000)(4320.000,2926.000)
\drawline(4200,2956)(4950,2956)
\drawline(2775,1531)(2776,795)
\drawline(2745.837,914.959)(2776.000,795.000)(2805.837,915.041)
\drawline(3945.839,3305.960)(3976.000,3186.000)(4005.839,3306.040)
\drawline(3976,3186)(3975,3931)
\drawline(1800,556)(2550,556)
\drawline(2430.000,526.000)(2550.000,556.000)(2430.000,586.000)
\drawline(1804.579,2026.102)(1741.000,1920.000)(1847.030,1983.700)
\drawline(1741,1920)(2616,2796)
\drawline(1735,1596)(2611,720)
\drawline(2504.934,783.640)(2611.000,720.000)(2547.360,826.066)
\drawline(375,781)(375,1531)
\drawline(405.000,1411.000)(375.000,1531.000)(345.000,1411.000)
\drawline(1350,556)(600,556)
\drawline(720.000,586.000)(600.000,556.000)(720.000,526.000)
\drawline(375,1981)(376,2723)
\drawline(405.838,2602.960)(376.000,2723.000)(345.838,2603.041)
\drawline(1416,1915)(544,2787)
\drawline(650.066,2723.360)(544.000,2787.000)(607.640,2680.934)
\drawline(1934.000,4185.000)(1814.000,4155.000)(1934.000,4125.000)
\drawline(1814,4155)(2560,4155)
\drawline(3750,4156)(3019,4155)
\drawline(3138.959,4185.164)(3019.000,4155.000)(3139.041,4125.164)
\drawline(2944,724)(3809,1589)
\drawline(3745.360,1482.934)(3809.000,1589.000)(3702.934,1525.360)
\drawline(544,3117)(1412,3987)
\drawline(1348.482,3880.861)(1412.000,3987.000)(1306.007,3923.238)
\drawline(4247.066,3928.360)(4141.000,3992.000)(4204.640,3885.934)
\thicklines
\drawline(4141,3992)(5016,3117)
\put(3075,1831){\makebox(0,0)[lb]{{\SetFigFont{6}{7.2}{\rmdefault}{\mddefault}{\updefault}$\begin{tabular}{r} 0\\-54\\-36\\-36 \end{tabular}\begin{tabular}{r} 0\\0\\-1\\1 \end{tabular}$}}}
\put(4200,1831){\makebox(0,0)[lb]{{\SetFigFont{6}{7.2}{\rmdefault}{\mddefault}{\updefault}$\begin{tabular}{r} 0\\-54\\-36\\0 \end{tabular}\begin{tabular}{r} -1\\0\\0\\2 \end{tabular}$}}}
\put(3450,856){\makebox(0,0)[lb]{{\SetFigFont{6}{7.2}{\rmdefault}{\mddefault}{\updefault}$\begin{tabular}{r} 0\\0\\0\\36 \end{tabular}\begin{tabular}{r} 1\\2\\-4\\0 \end{tabular}$}}}
\put(4800,2056){\makebox(0,0)[lb]{{\SetFigFont{6}{7.2}{\rmdefault}{\mddefault}{\updefault}$\begin{tabular}{r} 0\\-54\\-36\\0 \end{tabular}\begin{tabular}{r} 0\\2\\-3\\1 \end{tabular}$}}}
\put(3000,3031){\makebox(0,0)[lb]{{\SetFigFont{6}{7.2}{\rmdefault}{\mddefault}{\updefault}$\begin{tabular}{r} 72\\54\\54\\54 \end{tabular}\begin{tabular}{r} 0\\1\\-2\\1 \end{tabular}$}}}
\put(4275,2431){\makebox(0,0)[lb]{{\SetFigFont{6}{7.2}{\rmdefault}{\mddefault}{\updefault}$\begin{tabular}{r} 0\\-36\\-36\\0 \end{tabular}\begin{tabular}{r} 0\\1\\-1\\0 \end{tabular}$}}}
\put(2400,931){\makebox(0,0)[lb]{{\SetFigFont{6}{7.2}{\rmdefault}{\mddefault}{\updefault}$\begin{tabular}{r} -72\\0\\-24\\0 \end{tabular}\begin{tabular}{r} -1\\-2\\3\\1 \end{tabular}$}}}
\put(300,3706){\makebox(0,0)[lb]{{\SetFigFont{6}{7.2}{\rmdefault}{\mddefault}{\updefault}$\begin{tabular}{r} -72\\0\\-18\\-54 \end{tabular}\begin{tabular}{r} 0\\2\\-3\\1 \end{tabular}$}}}
\put(1800,31){\makebox(0,0)[lb]{{\SetFigFont{6}{7.2}{\rmdefault}{\mddefault}{\updefault}$\begin{tabular}{r} 72\\0\\36\\36 \end{tabular}\begin{tabular}{r} 0\\1\\-2\\2 \end{tabular}$}}}
\put(750,31){\makebox(0,0)[lb]{{\SetFigFont{6}{7.2}{\rmdefault}{\mddefault}{\updefault}$\begin{tabular}{r} -72\\0\\-36\\-36 \end{tabular}\begin{tabular}{r} -1\\-1\\2\\0 \end{tabular}$}}}
\put(1425,856){\makebox(0,0)[lb]{{\SetFigFont{6}{7.2}{\rmdefault}{\mddefault}{\updefault}$\begin{tabular}{r} 0\\0\\-12\\0 \end{tabular}\begin{tabular}{r} 1\\0\\0\\-2 \end{tabular}$}}}
\put(1950,4306){\makebox(0,0)[lb]{{\SetFigFont{6}{7.2}{\rmdefault}{\mddefault}{\updefault}$\begin{tabular}{r} 72\\0\\18\\54 \end{tabular}\begin{tabular}{r} -1\\-2\\4\\0 \end{tabular}$}}}
\put(3150,4306){\makebox(0,0)[lb]{{\SetFigFont{6}{7.2}{\rmdefault}{\mddefault}{\updefault}$\begin{tabular}{r} 72\\0\\18\\54 \end{tabular}\begin{tabular}{r} -1\\-1\\1\\1 \end{tabular}$}}}
\put(4650,3706){\makebox(0,0)[lb]{{\SetFigFont{6}{7.2}{\rmdefault}{\mddefault}{\updefault}$\begin{tabular}{r} 0\\-18\\0\\0 \end{tabular}\begin{tabular}{r} -1\\-1\\2\\1 \end{tabular}$}}}
\put(3600,3256){\makebox(0,0)[lb]{{\SetFigFont{6}{7.2}{\rmdefault}{\mddefault}{\updefault}$\begin{tabular}{r} -72\\-18\\-18\\-54 \end{tabular}\begin{tabular}{r} 1\\2\\-3\\-1 \end{tabular}$}}}
\put(2475,2056){\makebox(0,0)[lb]{{\SetFigFont{6}{7.2}{\rmdefault}{\mddefault}{\updefault}$\begin{tabular}{r} 72\\54\\60\\36 \end{tabular}\begin{tabular}{r} -1\\0\\0\\2 \end{tabular}$}}}
\put(0,2056){\makebox(0,0)[lb]{{\SetFigFont{6}{7.2}{\rmdefault}{\mddefault}{\updefault}$\begin{tabular}{r} -72\\0\\-36\\-36 \end{tabular}\begin{tabular}{r} 1\\1\\-3\\1 \end{tabular}$}}}
\put(0,856){\makebox(0,0)[lb]{{\SetFigFont{6}{7.2}{\rmdefault}{\mddefault}{\updefault}$\begin{tabular}{r} -72\\0\\-36\\-36 \end{tabular}\begin{tabular}{r} -1\\0\\0\\2 \end{tabular}$}}}
\put(1725,2431){\makebox(0,0)[lb]{{\SetFigFont{6}{7.2}{\rmdefault}{\mddefault}{\updefault}$\begin{tabular}{r} 0\\0\\6\\-18 \end{tabular}\begin{tabular}{r} -1\\-2\\3\\1 \end{tabular}$}}}
\put(900,2431){\makebox(0,0)[lb]{{\SetFigFont{6}{7.2}{\rmdefault}{\mddefault}{\updefault}$\begin{tabular}{r} 0\\0\\18\\-18 \end{tabular}\begin{tabular}{r} -1\\-1\\1\\1 \end{tabular}$}}}
\put(280,2896){\makebox(0,0)[lb]{{\SetFigFont{12}{14.4}{\rmdefault}{\mddefault}{\updefault}33}}}
\put(2700,4091){\makebox(0,0)[lb]{{\SetFigFont{12}{14.4}{\rmdefault}{\mddefault}{\updefault}57}}}
\put(5090,2896){\makebox(0,0)[lb]{{\SetFigFont{12}{14.4}{\rmdefault}{\mddefault}{\updefault}29}}}
\put(1485,1691){\makebox(0,0)[lb]{{\SetFigFont{12}{14.4}{\rmdefault}{\mddefault}{\updefault}17}}}
\put(2685,1696){\makebox(0,0)[lb]{{\SetFigFont{12}{14.4}{\rmdefault}{\mddefault}{\updefault}15}}}
\put(3890,1696){\makebox(0,0)[lb]{{\SetFigFont{12}{14.4}{\rmdefault}{\mddefault}{\updefault}26}}}
\put(5090,1691){\makebox(0,0)[lb]{{\SetFigFont{12}{14.4}{\rmdefault}{\mddefault}{\updefault}27}}}
\put(2695,501){\makebox(0,0)[lb]{{\SetFigFont{12}{14.4}{\rmdefault}{\mddefault}{\updefault}19}}}
\put(1535,496){\makebox(0,0)[lb]{{\SetFigFont{12}{14.4}{\rmdefault}{\mddefault}{\updefault}5}}}
\put(280,1696){\makebox(0,0)[lb]{{\SetFigFont{12}{14.4}{\rmdefault}{\mddefault}{\updefault}18}}}
\put(330,486){\makebox(0,0)[lb]{{\SetFigFont{12}{14.4}{\rmdefault}{\mddefault}{\updefault}6}}}
\put(3875,2896){\makebox(0,0)[lb]{{\SetFigFont{12}{14.4}{\rmdefault}{\mddefault}{\updefault}44}}}
\put(2675,2896){\makebox(0,0)[lb]{{\SetFigFont{12}{14.4}{\rmdefault}{\mddefault}{\updefault}16}}}
\put(3885,4096){\makebox(0,0)[lb]{{\SetFigFont{12}{14.4}{\rmdefault}{\mddefault}{\updefault}30}}}
\put(1495,4086){\makebox(0,0)[lb]{{\SetFigFont{12}{14.4}{\rmdefault}{\mddefault}{\updefault}58}}}
\end{picture}
}